\documentclass[12pt]{article}
\usepackage{latexsym}
\usepackage{amsmath}
\usepackage{amssymb}
\usepackage{amsfonts}
\usepackage{graphics,graphicx}
\usepackage[colorlinks=true]{hyperref}
\hypersetup{urlcolor=blue, citecolor=red}%
\usepackage[inline]{showlabels}
\usepackage{tikz}

\def\claim#1{\begin{trivlist}\item[\hskip\labelsep\bf#1]\it}
\def\endclaim{\end{trivlist}}

\numberwithin{equation}{section}

\newtheorem{theorem}{Theorem}[section]
\newtheorem{lemma}[theorem]{Lemma}
\newtheorem{pr}{Proposition}[section]

\newtheorem{definition}{Definition}[section]
\newtheorem{remark}[theorem]{Remark}

\newcommand{\eproof}{{\mbox{\ }~\hfill
\mbox{\large $\Box$} \par \vskip 10pt}}
\newcommand{\pf}{\noindent{\bf Proof}}
\newcommand{\R}{\mathbb R}
\newcommand{\supp}{\mbox{supp}\hspace{0.3mm}}

\title{Three-region inequalities for the second order elliptic equation with discontinuous coefficients and size estimate}
\author{E.~Francini\thanks{Universit\'a di Firenze, Italy. Email: francini@math.unifi.it}
\qquad C.-L.~Lin\thanks{National Cheng Kung University, Taiwan. Email: cllin2@mail.ncku.edu.tw}
\qquad S.~Vessella\thanks{Universit\'a di Firenze, Italy. Email: sergio.vessella@dmd.unifi.it}
\qquad J.-N.~Wang\thanks{National Taiwan University, Taiwan. Email: jnwang@ntu.edu.tw}}

\date{}

\begin{document}
\maketitle

\begin{abstract}
In this paper, we would like to derive a quantitative uniqueness estimate, the three-region inequality, for the second order elliptic equation with jump discontinuous coefficients. The derivation of the inequality relies on the Carleman estimate proved in our previous work \cite{dflvw15}. We then apply the three-region inequality to study the size estimate problem with one boundary measurement.  
\end{abstract}

\section{Introduction}

In this work we aim to study the size estimate problem with one measurement when the background conductivity has jump interfaces. A typical application of this study is to estimate the size of a cancerous tumor inside an organ by the electric impedance tomography (EIT).
In this case, considering discontinuous medium is typical, for instance, the conductivities of heart, liver, intestines are 0.70
(S/m), 0.10 (S/m), 0.03 (S/m), respectively. Previous works on this problem assumed that the conductivity of the studied body is Lipschitz continuous, see, for example, \cite{ar,ars}. The first result on the size estimate problem with a discontinuous background conductivity was given in \cite{nw13}, where only the two dimensional case was considered. In this paper, we will study the problem in dimension $n\ge 2$.  

The main ingredients of our method are quantitative uniqueness estimates for 
\begin{equation}\label{01}
\mbox{div}(A\nabla u)=0\quad\Omega\subset\R^n. 
\end{equation}
Those estimates are well-known when $A$ is Lipschitz continuous. The derivation of the estimates is based on the Carleman estimate or the frequency function method. For $n=2$ and $A\in L^\infty$, quantitative uniqueness estimates are obtained via the connection between \eqref{01} and quasiregular mappings. This is the method used in \cite{nw13}. For $n\ge 3$, the connection with quasiregular mappings is not true. Hence we return to the old method -- the Carleman estimate, to derive quantitative uniqueness estimates when $A$ is discontinuous. Precisely, when $A$ has a $C^{1,1}$ interface and is Lipschitz away from the interface, a Carleman estimate was obtained in \cite{dflvw15} (see \cite{ll, lr1, lr2} for related results). Here we will derive three-region inequalities using this Carleman estimate. The three-region inequality provides us a way to propagate "smallness" across the interface (see also \cite{lr1} for similar estimates). Relying on the three-region inequality, we then derive bounds of the size of an inclusion with one boundary measurement.  For other results on the size estimate, we mention \cite{amr} for the isotropic elasticity, \cite{mrv,mrv2,mrv3} for the isotropic/anisotropic thin plate, \cite{dlw,dlvw} for the shallow shell.

\section{The Carleman estimate}

In this section, we would like to describe the Carleman estimate derived in \cite{dflvw15}. We first denote $H_\pm=\chi_{\R^n_\pm}$ where $\R^n_\pm=\{(x,y)\in\R^{n-1}\times \R: y\gtrless0\}$ and $\chi_{\mathbb{R}^n_{\pm}}$ is the characteristic function of $\mathbb{R}^n_{\pm}$. Let $u_\pm\in C^\infty(\R^n)$ and define
\begin{equation*}\label{1.030}
u=H_+u_++H_-u_-=\sum_\pm H_{\pm}u_{\pm},
\end{equation*}
hereafter, $\sum_\pm a_\pm=a_++a_-$, and 
\begin{equation}\label{7.1}
\mathcal{L}(x,y,\partial)u:=\sum_{\pm}H_{\pm}{\rm div}_{x,y}(A_{\pm}(x,y)\nabla_{x,y}u_{\pm}),
\end{equation}
where
\begin{equation}\label{7.2}
A_{\pm}(x,y)=\{a^{\pm}_{ij}(x,y)\}^n_{i,j=1},\quad x\in \mathbb{R}^{n-1},y\in \mathbb{R}
\end{equation}
is a Lipschitz symmetric matrix-valued function satisfying, for given constants $\lambda_0\in (0,1]$, $M_0>0$,
\begin{equation}\label{7.3}
\lambda_0|z|^2\leq A_{\pm}(x,y)z\cdot z\leq \lambda^{-1}_0|z|^2,\, \forall (x,y)\in \mathbb{R}^n,\,\forall\, z\in \mathbb{R}^n
\end{equation}
and
\begin{equation}\label{7.4}
|A_{\pm}(x',y')-A_{\pm}(x,y)|\leq M_0(|x'-x|+|y'-y|).
\end{equation}
We write
\begin{equation}\label{7.5}
h_0(x):=u_+(x,0)-u_-(x,0),\ \forall\, x\in \mathbb{R}^{n-1},
\end{equation}
\begin{equation}\label{7.6}
h_1(x):=A_+(x,0)\nabla_{x,y}u_+(x,0)\cdot \nu-A_-(x,0)\nabla_{x,y}u_-(x,0)\cdot \nu,\ \forall\, x\in \mathbb{R}^{n-1},
\end{equation}
where $\nu=-e_n$.

For a function $h\in L^2(\mathbb{R}^{n})$, we define
\begin{equation*}
\hat{h}(\xi,y)=\int_{\mathbb{R}^{n-1}}h(x,y)e^{-ix\cdot\xi}\,dx,\quad \xi\in \mathbb{R}^{n-1}.
\end{equation*}
As usual $H^{1/2}(\mathbb{R}^{n-1})$ denotes the space of the functions $f\in L^2(\mathbb{R}^{n-1})$ satisfying
$$\int_{\mathbb{R}^{n-1}}|\xi||\hat{f}(\xi)|^2d\xi<\infty,$$
with the norm
\begin{equation}\label{semR}
||f||^2_{H^{1/2}(\mathbb{R}^{n-1})}=\int_{\mathbb{R}^{n-1}}(1+|\xi|^2)^{1/2}|\hat{f}(\xi)|^2d\xi.
\end{equation}
Moreover we define
$$[f]_{1/2,\mathbb{R}^{n-1}}=\left[\int_{\mathbb{R}^{n-1}}\int_{\mathbb{R}^{n-1}}\frac{|f(x)-f(y)|^2}{|x-y|^n}dydx\right]^{1/2},$$
and recall that there is a positive constant $C$, depending only on $n$, such that
\begin{equation*}
C^{-1}\int_{\mathbb{R}^{n-1}}|\xi||\hat{f}(\xi)|^2d\xi\leq[f]^2_{1/2,\mathbb{R}^{n-1}}\leq C\int_{\mathbb{R}^{n-1}}|\xi||\hat{f}(\xi)|^2d\xi,
\end{equation*}
so that the norm \eqref{semR} is equivalent to the norm $||f||_{L^2(\mathbb{R}^{n-1})}+[f]_{1/2,\mathbb{R}^{n-1}}$. From now on, we use the letters $C, C_0, C_1, \cdots$ to denote constants (depending on $\lambda_0,M_0,n$). The value of the constants may change from line to line, but it is always greater than $1$. We will denote by $B_r(x)$ the $(n-1)$-ball centered at $x\in \mathbb{R}^{n-1}$ with radius $r>0$. Whenever $x=0$ we denote $B_r=B_r(0)$.
\bigskip

\begin{theorem}\label{thm8.2}
Let $u$ and $A_{\pm}(x,y)$ satisfy \eqref{7.1}-\eqref{7.6}. There exist $L,\beta, \delta_0, r_0,\tau_0$ positive constants, with $r_0\leq 1$, depending on $\lambda_0, M_0, n$, such that if $\alpha_+>L\alpha_-$, $\delta\le\delta_0$ and $\tau\geq \tau_0$, then
\begin{equation}\label{8.24}
\begin{aligned}
&\sum_{\pm}\sum_{|k|=0}^2\tau^{3-2|k|}\int_{\mathbb{R}^n_{\pm}}|D^k{u}_{\pm}|^2e^{2\tau\phi_{\delta,\pm}(x,y)}dxdy+\sum_{\pm}\sum_{|k|=0}^1\tau^{3-2|k|}\int_{\mathbb{R}^{n-1}}|D^k{u}_{\pm}(x,0)|^2e^{2\phi_\delta(x,0)}dx\\
&+\sum_{\pm}\tau^2[e^{\tau\phi_\delta(\cdot,0)}u_{\pm}(\cdot,0)]^2_{1/2,\mathbb{R}^{n-1}}+\sum_{\pm}[D(e^{\tau\phi_{\delta,\pm}}u_{\pm})(\cdot,0)]^2_{1/2,\mathbb{R}^{n-1}}\\
\leq &C\left(\sum_{\pm}\int_{\mathbb{R}^n_{\pm}}|\mathcal{L}( x, y,\partial)(u_{\pm})|^2\,e^{2\tau\phi_{\delta,\pm}(x,y)}dxdy+[e^{\tau\phi_\delta(\cdot,0)}h_1]^2_{1/2,\mathbb{R}^{n-1}}\right.\\
&\left.+[\nabla_x(e^{\tau\phi_\delta}h_0)(\cdot,0)]^2_{1/2,\mathbb{R}^{n-1}}+\tau^{3}\int_{\mathbb{R}^{n-1}}|h_0|^2e^{2\tau\phi_\delta(x,0)}dx+\tau\int_{\mathbb{R}^{n-1}}|h_1|^2e^{2\tau\phi_\delta(x,0)}dx\right).
\end{aligned}
\end{equation}
where $u=H_+u_++H_-u_-$,  $u_{\pm}\in C^\infty(\mathbb{R}^{n})$ and ${\rm supp}\, u\subset B_{\delta/2}\times[-\delta r_0,\delta r_0]$, and $\phi_{\delta,\pm}(x,y)$ is given by
\begin{equation}\label{2.1}
\phi_{\delta,\pm}(x,y)=
\left\{
\begin{aligned}
\frac{\alpha_+y}{\delta}+\frac{\beta y^2}{2\delta^2}-\frac{|x|^2}{2\delta},\quad y\geq 0,\\
\frac{\alpha_-y}{\delta}+\frac{\beta y^2}{2\delta^2}-\frac{|x|^2}{2\delta},\quad y< 0,
\end{aligned}\right.
\end{equation}
and $\phi_\delta(x,0)=\phi_{\delta,+}(x,0)=\phi_{\delta,-}(x,0)$. 
\end{theorem}
\begin{remark}
It is clear that \eqref{8.24} remains valid if can add lower order terms\\ $\sum_{\pm}H_{\pm}\left(W\cdot\nabla_{x,y}u_{\pm}+Vu_{\pm}\right)$, where $W,V$ are bounded functions, to the operator ${\mathcal L}$. That is, one can substitute
\begin{equation}\label{8.242}
\mathcal{L}(x,y,\partial)u=\sum_{\pm}H_{\pm}{\rm div}_{x,y}(A_{\pm}(x,y)\nabla_{x,y}u_{\pm})+\sum_{\pm}H_{\pm}\left(W\cdot\nabla_{x,y}u_{\pm}+Vu_{\pm}\right)
\end{equation}
in \eqref{8.24}.
\end{remark}

\section{Three-region inequalities}

Based on the Carleman estimate given in Theorem~\ref{thm8.2}, we will derive three-region inequalities across the interface $y=0$. Here we consider $u=H_+u_++H_-u_-$ satisfying
\[
\mathcal{L}(x,y,\partial)u=0\quad\mbox{in}\quad\R^n,
\]
where ${\mathcal L}$ is given in \eqref{8.242} and
\[
\|W\|_{L^\infty(\R^n)}+\|V\|_{L^\infty(\R^n)}\le\lambda_0^{-1}.
\]
Fix any $\delta\le\delta_0$, where $\delta_0$ is given in Theorem~\ref{thm8.2}.

\begin{theorem}\label{thm9.1}
Let $u$ and $A_{\pm}(x,y)$ satisfy \eqref{7.1}-\eqref{7.6} with $h_0=h_1=0$. Then there exist $C$ and ${R}$, depending only on $\lambda_0, M_0, n$, such that if $\ 0<R_1,R_2\leq R$, then
\begin{equation}\label{9.1}
\int_{U_2}|u|^2 dx \le (e^{\tau_0R_2}+CR_1^{-4}) \left(\int_{U_1}|u|^2dxdy\right)^{\frac{R_2}{2R_1+3R_2}}\left(\int_{U_3}|u|^2dxdy\right)^{\frac{2R_1+2R_2}{2R_1+3R_2}},
\end{equation}
where $\tau_0$ is the constant derived in Theorem~\ref{thm8.2},
\[
\begin{aligned}
&U_1=\{z\geq-4R_2,\,\frac{R_1}{8a}<y<\frac{R_1}{a}\},\\
&U_2=\{-R_2\leq z\leq\frac{R_1}{2a},\,y<\frac{R_1}{8a}\},\\
&U_3=\{z\geq-4R_2,\,y<\frac{R_1}{a}\},
\end{aligned}
\]
$a=\alpha_+/\delta$, and
\begin{equation}\label{zxy}
z(x,y)=\frac{\alpha_-y}{\delta}+\frac{\beta y^2}{2\delta^2}-\frac{|x|^2}{2\delta}.
\end{equation}
\end{theorem}
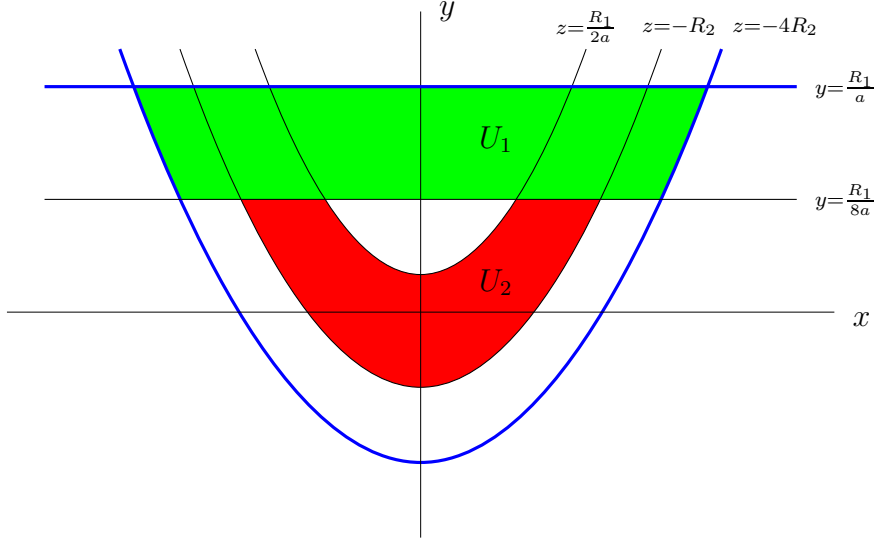
\begin{figure}[ht]
\centering
\begin{tikzpicture}
\begin{scope}
\clip (-5,3) rectangle (5,1.5);
\fill[green] (-4,3.5) parabola bend (0,-2) (4,3.5);
\end{scope}
\draw[line width=1.2pt, blue] (-4,3.5) parabola bend (0,-2) (4,3.5);
\begin{scope}
\clip (-5,1.5) rectangle (5,-1.2);
\fill[red] (-3.2,3.5) parabola bend (0,-1) (3.2,3.5);
\fill[white] (-2.2,3.5) parabola bend (0,0.5) (2.2,3.5);
\end{scope}
\draw[line width=0.2pt]  (5.5,0) -- (-5.5,0);
\node[right] at (0.1,4) {$y$}; 
\draw[line width=0.2pt]  (0,4) -- (0,-3);
\node[right] at (5.6,-0.1) {$x$};
\node[right] at (4.0,3.8) {$\scriptstyle z=-4R_2$}; 
\draw (-3.2,3.5) parabola bend (0,-1) (3.2,3.5);
\node[right] at (2.8,3.8) {$\scriptstyle z=-R_2$};
\draw (-2.2,3.5) parabola bend (0,0.5) (2.2,3.5);
\draw[line width=1.2pt, blue] (-5,3) -- (5,3);
\draw (-5,1.5) -- (5,1.5);
\node at (2.2,3.8) {$\scriptstyle z=\frac{R_1}{2a}$};
\node[right] at (5.1,3) {$\scriptstyle y=\frac{R_1}{a}$};
\node[right] at (5.1,1.5) {$\scriptstyle y=\frac{R_1}{8a}$};
\node at (1,2.3) {$U_1$};
\node at (1,0.4) {$U_2$};
\end{tikzpicture}
\caption{$U_1$ and $U_2$ are shown in green and red, respectively. $U_3$ is the region enclosed by blue boundaries. }
\label{fig1}
\end{figure}

\pf. To apply the estimate \eqref{8.24}, $u$ needs to satisfy the support condition. Also, we can choose $\alpha_+$ and $\alpha_-$ in Theorem~\ref{thm8.2} such that $\alpha_+>\alpha_-$. We can choose $r\le r_0$ satisfying
\begin{equation}\label{1.5}
r\le\min\left\{\frac{13\alpha_-}{8\beta},\frac{2\delta}{19\alpha_-+8\beta}\right\}.
\end{equation}
Note that the choices of $\delta, r$ also depend on $\lambda_0,M_0,n$. We then set
\[
R=\frac{\alpha_- r}{16}.
\]
It follows from \eqref{1.5} that
\begin{equation}\label{bee}
R\le\frac{13\alpha_-^2}{128\beta}.
\end{equation}
Given $0<R_1<R_2\le R$. Let $\vartheta_1(t)\in C^{\infty}_0 ({\mathbb R})$ satisfy $0\le\vartheta_1(t)\leq 1$ and
\begin{equation*}
\vartheta_1 (t)=\left\{
\begin{aligned}
&1,\quad t>-2R_2,\\
&0,\quad t\leq -3R_2.
\end{aligned}\right.
\end{equation*}
Also, define $\vartheta_2(y)\in C^{\infty}_0 ({\mathbb R})$ satisfying $0\le\vartheta_2(y)\leq 1$ and
\begin{equation*}
\vartheta_2 (y)=\left\{
\begin{aligned}
&0,\quad y\geq \frac{R_1}{2a},\\
&1,\quad y<\frac{R_1}{4a}.
\end{aligned}\right.
\end{equation*}
Finally, we define $\vartheta(x,y)=\vartheta_1(z(x,y))\vartheta_2(y)$, where $z$ is defined by \eqref{zxy}. 

We now check the support condition for $\vartheta$. From its definition, we can see that $\supp\vartheta$ is contained in
\begin{equation}\label{2.2}
\left\{
\begin{aligned}
&z(x,y)=\frac{\alpha_-y}{\delta}+\frac{\beta y^2}{2\delta^2}-\frac{|x|^2}{2\delta}>-3R_2,\\
& y<\frac{R_1}{2a}.
\end{aligned}\right.
\end{equation}
In view of the relation
\[
\alpha_+>\alpha_-\quad\mbox{and}\quad a=\frac{\alpha_+}{\delta},
\]
we have that
\begin{equation*}
\frac{R_1}{2a}<\frac{\delta}{2\alpha_-}\cdot R_1<\frac{\delta}{\alpha_-}\cdot\frac{\alpha_-r}{16}<\delta r,
\end{equation*}
i.e., $y<\delta r\le\delta r_0$. Next, we  observe that
\[
-3R_2>-3R=-\frac{3\alpha_-r}{16}>\frac{\alpha_-}{\delta}(-\delta r)+\frac{\beta}{2\delta^2}(-\delta r)^2,
\]
which gives $-\delta r<y$ due to  \eqref{1.5}. Consequently, we verify that $|y|<\delta r$. One the other hand, from the first condition of \eqref{2.2} and \eqref{1.5}, we see that
\begin{equation*}
\frac{|x|^2}{2\delta}<3R_2+\frac{\alpha_-y}{\delta}+\frac{\beta y^2}{2\delta^2}\leq\frac{3\alpha_-r}{16}+\frac{\alpha_-}{\delta}\cdot \delta r+\frac{\beta}{2\delta^2}\cdot \delta^2r^2\leq\frac{\delta}{8},
\end{equation*}
which gives $|x|<\delta/2$.

Since $h_0=0$, we have that
\begin{equation}\label{8.90}
\vartheta(x,0)u_+(x,0)-\vartheta(x,0)u_-(x,0)=0,\;\forall\; x\in\R^{n-1}.
\end{equation}
Applying \eqref{8.24} to $\vartheta u$ and using \eqref{8.90} yields
\begin{equation}\label{9.2}
\begin{aligned}
&\sum_{\pm}\sum_{|k|=0}^2\tau^{3-2|k|}\int_{\mathbf{R}^n_{\pm}}|D^k(\vartheta u_{\pm})|^2e^{2\tau\phi_{\delta,\pm}(x,y)}dxdy\\
\le &C\sum_{\pm}\int_{\mathbf{R}^n_{\pm}}|{\mathcal L}( x, y,\partial)(\vartheta u_{\pm})|^2\,e^{2\tau\phi_{\delta,\pm}(x,y)}dxdy\\
&+C\tau\int_{\mathbf{R}^{n-1}}|A_+(x,0)\nabla_{x,y}(\vartheta u_+(x,0))\cdot \nu-A_-(x,0)\nabla_{x,y}(\vartheta u_-)(x,0)\cdot \nu|^2e^{2\tau\phi_\delta(x,0)}dx\\
&+C[e^{\tau\phi_\delta(\cdot,0)}\big(A_+(x,0)\nabla_{x,y}(\vartheta u_+)(x,0)\cdot \nu-A_-(x,0)\nabla_{x,y}(\vartheta u_-)(x,0)\cdot \nu\big)]^2_{1/2,\mathbf{R}^{n-1}}.
\end{aligned}
\end{equation}
We now observe that $\nabla_{x,y}\vartheta_1(z)=\vartheta_1'(z)\nabla_{x,y}z=\vartheta_1'(z)(-\frac{x}{\delta},\frac{\alpha_-}{\delta}+\frac{\beta y}{\delta^2})$ and it is nonzero only when 
\[
-3R_2<z<-2R_2.
\]
Therefore, when $y=0$, we have
\[
2R_2<\frac{|x|^2}{2\delta}<3R_2.
\]
Thus, we can see that
\begin{equation}\label{1001}
|\nabla_{x,y}\vartheta(x,0)|^2\le CR_2^{-2}\left(\frac{6R_2}{\delta}+\frac{\alpha_-^2}{\delta^2}\right)\le C R_2^{-2}.
\end{equation}
By $h_0(x)=h_1(x)=0$, \eqref{1001}, and the easy estimate of \cite[Proposition~4.2]{dflvw15}, it is not hard to estimate
\begin{equation}\label{9.3}
\begin{aligned}
&\tau\int_{\mathbf{R}^{n-1}}|A_+(x,0)\nabla_{x,y}(\vartheta u_+(x,0))\cdot \nu-A_-(x,0)\nabla_{x,y}(\vartheta u_-)(x,0)\cdot \nu|^2e^{2\tau\phi_\delta(x,0)}dx\\
&+[e^{\tau\phi_\delta(\cdot,0)}\big(A_+(x,0)\nabla_{x,y}(\vartheta u_+)(x,0)\cdot \nu-A_-(x,0)\nabla_{x,y}(\vartheta u_-)(x,0)\cdot \nu\big)]^2_{1/2,\mathbf{R}^{n-1}}\\
\le &\,C R_2^{-2}e^{-4\tau R_2}\left(\tau\int_{\{\sqrt{4\delta R_2}\leq|x|\leq\sqrt{6\delta R_2}\}}|u_+(x,0)|^2dx+[u_+(x,0)]^2_{1/2,\{\sqrt{4\delta R_2}\leq|x|\leq\sqrt{6\delta R_2}\}}\right)\\
&+C\tau^2R_2^{-3}e^{-4\tau R_2}\int_{\{\sqrt{4\delta R_2}\leq|x|\leq\sqrt{6\delta R_2}\}}|u_+(x,0)|^2dx\\
\le &\,C\tau^2R_2^{-3}e^{-4\tau R_2}E,
\end{aligned}
\end{equation}
where 
\[
E=\int_{\{\sqrt{4\delta R_2}\leq|x|\leq\sqrt{6\delta R_2}\}}|u_+(x,0)|^2dx+[u_+(x,0)]^2_{1/2,\{\sqrt{4\delta R_2}\leq|x|\leq\sqrt{6\delta R_2}\}}.
\]

Expanding ${\cal L}( x, y,\partial)(\vartheta u_{\pm})$ and considering the set where $D\vartheta\neq0$, we can estimate
\begin{equation}\label{9.4}
\begin{aligned}
&\sum_{\pm}\sum_{|k|=0}^1\tau^{3-2|k|}\int_{\{-2R_2\leq z\leq \frac{R_1}{2a},\,y<\frac{R_1}{4a}\}}|D^ku_{\pm}|^2e^{2\tau\phi_{\delta,\pm}(x,y)}dxdy\\
\le &\,C\sum_{\pm}\sum_{|k|=0}^1R_2^{2(|k|-2)}\int_{\{-3R_2\leq z\leq-2R_2,\,y<\frac{R_1}{2a}\}}|D^ku_{\pm}|^2e^{2\tau\phi_{\delta,\pm}(x,y)}dxdy\\
&\,+C\sum_{|k|=0}^1R_1^{2(|k|-2)}\int_{\{-3R_2\leq z,\,\frac{R_1}{4a}< y<\frac{R_1}{2a}\}}|D^ku_{+}|^2e^{2\tau\phi_{\delta,+}(x,y)}dxdy\\
&\,+C\tau^2R_2^{-3}e^{-4\tau R_2}E\\
\le &\, C\sum_{\pm}\sum_{|k|=0}^1R_2^{2(|k|-2)}e^{-4\tau R_2}e^{2\tau\frac{(\alpha_+-\alpha_-)}{\delta}\frac{R_1}{4a}}\int_{\{-3R_2\leq z\leq-2R_2,\,y<\frac{R_1}{4a}\}}|D^{k}u_\pm|^2dxdy\\
&+\sum_{|k|=0}^1R_1^{2(|k|-2)}e^{2\tau\frac{\alpha_+}{\delta}\frac{R_1}{2a}}e^{2\tau\frac{\beta}{2\delta^2}(\frac{R_1}{2a})^2}\int_{\{z\geq-3R_2,\,\frac{R_1}{4a}<y<\frac{R_1}{2a}\}}|D^{k}u_+|^2dxdy\\
&\,+C\tau^2R_2^{-3}e^{-4\tau R_2}E.
\end{aligned}
\end{equation}
Let us denote $U_1=\{z\geq-4R_2,\,\frac{R_1}{8a}<y<\frac{R_1}{a}\}$,
$U_2=\{-R_2\leq z\leq\frac{R_1}{2a},\,y<\frac{R_1}{8a}\}$. From \eqref{9.4} and interior estimates (Caccioppoli's type inequality), we can derive that
\begin{equation}\label{9.5}
\begin{aligned}
&\tau^3e^{-2\tau R_2}\int_{U_2}|u|^2dxdy\\
\leq&\,\tau^{3}e^{-2\tau R_2}\int_{\{-R_2\leq z\leq\frac{R_1}{2a},\,y<\frac{R_1}{8a}\}}|u|^2dxdy\\
\leq&\,\sum_{\pm}\tau^{3}\int_{\{-2R_2\leq z\leq\frac{R_1}{2a},\,y<\frac{R_1}{4a}\}}| u_{\pm}|^2e^{2\tau\phi_{\delta,\pm}(x,y)}dxdy\\
\le &\, C\sum_{\pm}\sum_{|k|=0}^1R_2^{2(|k|-2)}e^{-4\tau R_2}e^{2\tau\frac{(\alpha_+-\alpha_-)}{\delta}\frac{R_1}{4a}}\int_{\{-3R_2\leq z\leq-2R_2,\,y<\frac{R_1}{4a}\}}|D^{k}u_\pm|^2dxdy\\
&+\sum_{|k|=0}^1R_1^{2(|k|-2)}e^{2\tau\frac{\alpha_+}{\delta}\frac{R_1}{2a}}e^{2\tau\frac{\beta}{2\delta^2}(\frac{R_1}{2a})^2}\int_{\{z\geq-3R_2,\,\frac{R_1}{4a}<y<\frac{R_1}{2a}\}}|D^{k}u_+|^2dxdy\\
&\,+C\tau^2R_2^{-3}e^{-4\tau R_2}E\\
\le &\,CR_1^{-4}e^{-3\tau R_2}\int_{\{-4R_2\leq z\leq-R_2,\,y<\frac{R_1}{a}\}}| u|^2dxdy+C\tau^2R_2^{-3}e^{-4\tau R_2}E\\
&\,+CR_1^{-4}e^{(1+\frac{\beta R_1}{4\alpha_-^2})\tau R_1}\int_{\{z\geq-4R_2,\,\frac{R_1}{8a}<y<\frac{R_1}{a}\}}| u|^2dxdy\\
\le &CR_1^{-4}\left(e^{2\tau R_1}\int_{U_1}|u|^2dxdy+\tau^2e^{-3\tau R_2}F\right),
\end{aligned}
\end{equation}
where \[F=\int_{\{z\geq-4R_2,\,y<\frac{R_1}{a}\}}| u|^2dxdy\] and we used the inequality $\frac{\beta R_1}{4\alpha_-^2}<1$ due to \eqref{bee}.

Dividing $\tau^3e^{-2\tau R_2}$ on both sides of \eqref{9.5} implies that
\begin{equation}\label{9.6}
\int_{U_2}|u|^2dxdy\le CR_1^{-4}\left( e^{2\tau (R_1+R_2)}\int_{U_1}|u|^2dxdy+e^{-\tau R_2}F\right).
\end{equation}
Now, we consider two cases. If $\int_{U_1} |u|^2dxdy\ne 0$ and
$$e^{2\tau_0 (R_1+R_2)}\int_{U_1}|u|^2dxdy<e^{-\tau_0 R_2}F,$$
then we can pick a $\tau>\tau_0$ such that
$$
e^{2\tau (R_1+R_2)}\int_{U_1}|u|^2dxdy=e^{-\tau R_2}F.
$$
Using such $\tau$, we obtain from \eqref{9.6} that
\begin{equation}\label{9.7}
\begin{aligned}
\int_{U_2}|u|^2dxdy&\le CR_1^{-4} e^{2\tau (R_1+R_2)}\int_{U_1}|u|^2dxdy\\
&=CR_1^{-4}\left(\int_{U_1}|u|^2dxdy\right)^{\frac{R_2}{2R_1+3R_2}}(F)^{\frac{2R_1+2R_2}{2R_1+3R_2}}.
\end{aligned}
\end{equation}
If $\int_{U_1}|u|^2dxdy= 0$, then letting $\tau\to\infty$ in \eqref{9.6} we have $\int_{U_2}|u|^2dxdy=0$ as well. The
three-regions inequality \eqref{9.1} obviously holds.

On the other hand, if
$$ e^{-\tau_0 R_2}F\leq e^{2\tau_0 (R_1+R_2)}\int_{U_1}|u|^2dxdy,$$
then we have
\begin{equation}\label{9.8}
\begin{aligned}
\int_{U_2}|u|^2 dx&\leq \left(F\right)^{\frac{R_2}{2R_1+3R_2}}\left(F\right)^{\frac{2R_1+2R_2}{2R_1+3R_2}}\\
&\leq \exp{(\tau_0R_2)}\left(\int_{U_1}|u|^2dxdy\right)^{\frac{R_2}{2R_1+3R_2}}\left(F\right)^{\frac{2R_1+2R_2}{2R_1+3R_2}}.
\end{aligned}
\end{equation}
Putting together \eqref{9.7}, \eqref{9.8}, we arrive at
\begin{equation}\label{9.9}
\int_{U_2}|u|^2 dx \le
(\exp{(\tau_0R_2)}+CR_1^{-4})\left(\int_{U_1}|u|^2dxdy\right)^{\frac{R_2}{2R_1+3R_2}}\left(F\right)^{\frac{2R_1+2R_2}{2R_1+3R_2}}.
\end{equation}
\eproof

\section{Size estimate}

We will apply the three-region inequality \eqref{9.1} to estimate the size of embedded inclusion in this section. Here we denote $\Omega$ a bounded open set in $\R^n$ with $C^{1,\alpha}$ boundary $\partial\Omega$ with constants $s_0,L_0$, where $0<\alpha\le 1$. Assume that $\Sigma$ is a $C^{2}$ hypersurface with constants $r_0,K_0$ satisfying
\begin{equation}\label{d0}
\mbox{dist}(\Sigma,\partial\Omega)\ge d_0
\end{equation}
for some $d_0>0$. We divide $\Omega$ into three sets, namely, 
\[
\Omega=\Omega_+\cup\Sigma\cup \Omega_-
\] 
where $\Omega_\pm$ are open subsets. Note that $\overline\Omega_-=\partial\Omega\cup\Sigma$ and $\partial\Omega_+=\Sigma$. We also define
\[
\Omega_h=\{x\in\Omega:\mbox{dist}(x,\partial\Omega)>h\}.
\]
\begin{definition}{\rm [}$C^{1,\alpha}$ \rm{regularity}{]}
We say that $\Sigma$ is $C^{2}$ with constants $r_0, K_0$ if for any $P\in\Sigma$ there exists a rigid transformation of coordinates under which $P=0$ and
\[
\Omega_\pm\cap B(0,r_0)=\{(x,y)\in B(0,r_0)\subset\R^n: y\gtrless\psi(x)\},
\]
where $\psi$ is a $C^{2}$ function on $B_{r_0}(0)$ satisfying $\psi(0)=0$ and 
\[
\|\psi\|_{C^{2}(B_{r_0}(0))}\le K_0.
\]
The definition of $C^{1,\alpha}$ boundary is similar. Note that $B(a,r)$ stands for the $n$-ball centered at $a$ with radius $r>0$. We remind the reader that $B_r(a)$ denotes the $(n-1)$-ball centered at $a$ with radius $r>0$. 
\end{definition}

Assume that $A_{\pm}=\{a_{ij}^\pm(x,y)\}_{i,j=1}^n$ satisfy \eqref{7.3} and \eqref{7.4}. Let us define $H_\pm=\chi_{\Omega_\pm}$, $A=H_+A_++H_-A_-$, $u=H_+u_++H_-u_-$. We now consider the conductivity equation
\begin{equation}\label{aeq}
\mbox{div}(A\nabla u)=0\quad\mbox{in}\quad\Omega.
\end{equation}
It is not hard to check that $u$ satisfies homogeneous transmission conditions \eqref{7.5}, \eqref{7.6} (with $h_0=h_1=0$), where in this case $\nu$ is the outer normal of $\Sigma$. For $\phi\in H^{1/2}(\partial\Omega)$, let $u$ solve \eqref{aeq} and satisfy the boundary value $u=\phi$ on $\partial\Omega$.

Next we assume that $D$ is a measurable subset of $\Omega$. Suppose that $\hat A$ is a symmetric $n\times n$ matrix with $L^\infty(\Omega)$ entries.  In addition, we assume that there exist $\eta>0, \zeta>1$ such that
\begin{equation}\label{jump1}
(1+\eta)A\le\hat A\le\zeta A\quad\mbox{a.e. in}\quad\Omega
\end{equation} 
or $\eta>0, 0<\zeta<1$ such that
\begin{equation}\label{jump2}
\zeta A\le\hat A\le(1-\eta) A\quad\mbox{a.e. in}\quad\Omega.
\end{equation} 
Now let $v=H_+v_++H_-v_-$ be the solution of
\begin{equation}\label{bvp2}
\left\{
\begin{aligned}
&\mbox{div}((A\chi_{\Omega\setminus\bar D}+\hat A\chi_D)\nabla v)=0\quad\mbox{in}\quad\Omega,\\
&v=\phi\quad\mbox{on}\quad\partial\Omega.
\end{aligned}\right.
\end{equation}
The inverse problem considered here is to estimate $|D|$ by the knowledge of $\{\phi, A\nabla v\cdot\nu|_{\partial\Omega}\}$. In this work we would like to consider the most interesting case where 
\begin{equation}\label{interior}\bar{D}\subseteq\bar{\Omega}_+.\end{equation} 
In practice, one could think of $\Omega_+$ being an organ and $D$ being a tumor. The aim is to estimate the size of $D$ by measuring one pair of voltage and current on the surface of the body. 

We denote $W_0$ and $W$ the powers required to maintain the voltage $\phi$ on $\partial\Omega$ when the inclusion $D$ is absent or present. It is easy to see that
\[
W_0=\int_{\partial\Omega}\phi A\nabla u\cdot\nu=\int_{\Omega}A\nabla u\cdot\nabla u
\]
and 
\[
W=\int_{\partial\Omega}\phi (A\chi_{\Omega\setminus\bar D}+\hat A\chi_D)\nabla v\cdot\nu=\int_{\Omega}(A\chi_{\Omega\setminus\bar D}+\hat A\chi_D)\nabla v\cdot\nabla v.
\]

The size of $D$ will be estimate by the power gap $W-W_0$. To begin, we recall the following energy inequalities proved in \cite{ars}.
\begin{lemma}\cite[Lemma~2.1]{ars}\label{energy} Assume that $A$
satisfies the ellipticity condition \eqref{7.3}. If either \eqref{jump1}
or \eqref{jump2} holds, then
\begin{equation}
C_{1}\int_{D}|\nabla u|^{2}\le\left|W_{0}-W\right|\le C_{2}\int_{D}|\nabla u|^{2},\label{ineq1}
\end{equation}
where $C_{1},C_{2}$ are constants depending only on $\lambda$, $\eta$,
and $\zeta$. \end{lemma}
The derivation of bounds on $|D|$ will be based on \eqref{ineq1} and the following Lipschitz propagation of smallness for $u$. 
\begin{pr}(Lipschitz propagation of smallness)\label{lippro}
Let $u\in H^{1}(\Omega)$ be the solution of \eqref{aeq} with Dirichlet data $\phi$.   For any $B(x,\rho)\subset\Omega_+$, we have that
\begin{equation}
\int_{B(x,\rho)}|\nabla u|^{2}\ge C\int_{\Omega}|\nabla u|^{2},\label{lp}
\end{equation}
where $C$ depends on $\Omega_\pm$, $d_0$, $\lambda_0$, $M_{0}$, $r_0$, $K_0$, $s_0$, $L_0$, $\alpha$, $\alpha'$, $\rho$, and 
$$\frac{\|\phi-\phi_0\|_{C^{1,\alpha'}(\partial\Omega)}}{\|\phi-\phi_0\|_{H^{1/2}(\partial\Omega)}},$$
for $\phi_0=|\partial\Omega|^{-1}\int_{\partial\Omega}\phi$. Here $\alpha'$ satisfies $0<\alpha'<\frac{\alpha}{(\alpha+1)n}$. 
\end{pr}

Before proving Proposition~\ref{lippro}, we need to adjust the three-region inequality \eqref{9.1} for the $C^{2}$ interface $\Sigma$. Let $0\in\Sigma$ and the coordinate transform $(x',y')=T(x,y)=(x,y-\psi(x))$ for $x\in B_{s_0}(0)$. Denote $\tilde U=T(B(0,s_0))$ and $\tilde{\mathcal{A}}_\pm=\{\tilde a_{i,j}^{\pm}\}_{i,j=1}^n$ the coefficients of $A_\pm$ in the new coordinates $(x',y')$. It is easy to see that 
$\tilde{\mathcal{A}}_\pm$ satisfies \eqref{7.3} and \eqref{7.4} with possible different constants $\tilde\lambda_0, \tilde M_0$, depending on $\lambda_0, M_0, r_0, K_0$. Then there exist $C$ and $\tilde R$, depending on $\tilde\lambda_0, \tilde M_0, n$, such that for 
\begin{equation}\label{r1r2}
0<R_1<R_2\le\tilde R
\end{equation}
and $U_1, U_2, U_3$ defined as in Theorem~\ref{thm9.1}, we have that $U_3\subset\tilde U$ (so $U_1, U_2$ are contained in $\tilde U$ as well) and \eqref{9.1} holds. Now let $\tilde U_j=T^{-1}(U_j)$, $j=1,2,3$, then \eqref{9.1} becomes
\begin{equation}\label{3r}
\int_{\tilde U_2}|u|^2 dxdy \le C \left(\int_{\tilde U_1}|u|^2dxdy\right)^{\frac{R_2}{2R_1+3R_2}}\left(\int_{\tilde U_3}|u|^2dxdy\right)^{\frac{2R_1+2R_2}{2R_1+3R_2}},
\end{equation}
where $C$ depends on $\lambda_0, M_0, r_0, K_0, n, R_1, R_2$. Furthermore, by Caccioppoli's inequality and generalized Poincar\'e's inequality (see (3.8) in \cite{amr08}), we obtain from \eqref{3r} that
\begin{equation}\label{3rd}
\int_{\tilde U_2}|\nabla u|^2 dxdy \le C \left(\int_{\tilde U_1}|\nabla u|^2dxdy\right)^{\frac{R_2}{2R_1+3R_2}}\left(\int_{\tilde U_3}|\nabla u|^2dxdy\right)^{\frac{2R_1+2R_2}{2R_1+3R_2}}
\end{equation}
with a possibly different constant $C$. 

Since $A_+$ (respectively $A_-$) is Lipschitz in $\Omega_+$ (respectively $\Omega_-$), the following three-sphere inequality is well-known. Let $u_\pm$ be a solution to $\mbox{div}(A_\pm\nabla u_\pm)=0$ in $\Omega_\pm$. Then for $B(x_0,\bar r)\subset\Omega_+$ (or $B(x_0,\bar r)\subset\Omega_-$) and $0<r_1<r_2<r_3<\bar r$, we have that
\begin{equation}\label{3s}
\int_{B(x_0,r_2)}|\nabla u_\pm|^2dxdy\le C\left(\int_{B(x_0,r_1)}|\nabla u_\pm|^2dxdy\right)^{\theta}\left(\int_{B(x_0,r_3)}|\nabla u_\pm|^2dxdy\right)^{1-\theta},
\end{equation} 
where $0<\theta<1$ and $C$ depend on $\lambda_0,M_0,n, r_1/r_3, r_2/r_3$. 

Now we are ready to prove Proposition~\ref{lippro}.

\smallskip
\noindent{\bf Proof of Proposition~\ref{lippro}}. It suffices to study the case where $\rho$ is small. Since $\Sigma\in C^2$, it satisfies both the uniform interior and exterior sphere properties, i.e., there exists $a_0>0$ such that for all $z\in\Sigma$, there exist balls $B\subset\Omega_+$ and $B'\subset\Omega_-$ of radius $a_0$ such that $\overline B\cap\Sigma=\overline B'\cap\Sigma=\{z\}$. Next let $\nu_z$ be the unit normal at $z\in\Sigma$ pointing into $\Omega_+$ (inwards) and $L=\{z+t\nu_z\subset\R^n: t\in[\rho_0,-3\rho_0]\}$. We then fix $R_1,R_2$ satisfying \eqref{r1r2} and choose $\rho_0>0$ so that
\[
S_z=\cup_{y\in L}B(y,\rho_0)\subset\tilde U_2.
\] 
Denote $\kappa=R_2/(2R_1+3R_2)$. Note that we move the construction of the three-region inequality from $0$ to $z$.

Let $x\in\Omega_+$ and consider $B(x,\rho)\subset\Omega_+$, where $\rho\le\min\{a_0,\rho_0\}$. For any $y\in\Omega_{2\rho}$, we discuss three cases.\\
(i) Let $y\in\Omega_{+,\rho}$, then by \eqref{3s} and the chain of balls argument, we have that
\begin{equation}
\frac{\int_{B(y,\rho)}|\nabla u|^2}{\int_\Omega|\nabla u|^2}\le C\left(\frac{\int_{B(x,\rho)}|\nabla u|^2}{\int_\Omega|\nabla u|^2}\right)^{\theta^{N_1}},\label{e1}
\end{equation}
where $N_1$ depends on $\Omega_+$ and $\rho$. \\
(ii) Let $y\in\{y\in\overline\Omega_+:\mbox{dist}(y,\Sigma)\le\rho\}\cup\{y\in\Omega_-:\mbox{dist}(y,\Sigma)\le 3\rho\}$, then $B(y,\rho)\subset S_z$ for some $z\in\Sigma$. Note that $\tilde U_1\subset\Omega_{+,\rho}$ (taking $\rho$ even smaller if necessary). We then apply \eqref{e1} iteratively to estimate
\begin{equation}\label{e2}
\frac{\int_{\tilde U_1}|\nabla u|^2}{\int_\Omega|\nabla u|^2}\le C\left(\frac{\int_{B(x,\rho)}|\nabla u|^2}{\int_\Omega|\nabla u|^2}\right)^{\theta^{N_1}},
\end{equation}
where $C$ depends on $\tilde U_1$ and $\rho$. Combining estimates \eqref{e2} and \eqref{3rd} yields
\begin{equation}\label{e3}
\frac{\int_{B(y,\rho)}|\nabla u|^2}{\int_\Omega|\nabla u|^2}\le C\left(\frac{\int_{B(x,\rho)}|\nabla u|^2}{\int_\Omega|\nabla u|^2}\right)^{\kappa\theta^{N_1}}.
\end{equation} 
\\
(iii) Finally, we consider the case where $y\in\Omega_-\cap\Omega_{2\rho}$ and $\mbox{dist}(y,\Sigma)>3\rho$. We observe that if $y_\ast=z+(-3\rho)\nu_z$, then \eqref{e3} implies
\begin{equation}\label{e33}
\frac{\int_{B(y_\ast,\rho)}|\nabla u|^2}{\int_\Omega|\nabla u|^2}\le C\left(\frac{\int_{B(x,\rho)}|\nabla u|^2}{\int_\Omega|\nabla u|^2}\right)^{\kappa\theta^{N_1}}.
\end{equation}
Again using \eqref{3s} and the chain of balls argument (starting with \eqref{e33}), we obtain that
\begin{equation}\label{e5}
\frac{\int_{B(y,\rho)}|\nabla u|^2}{\int_\Omega|\nabla u|^2}\le C\left(\frac{\int_{B(x,\rho)}|\nabla u|^2}{\int_\Omega|\nabla u|^2}\right)^{\kappa\theta^{N_1}\theta^{N_2}}.
\end{equation}
Putting together \eqref{e1}, \eqref{e3}, and \eqref{e5} gives
\begin{equation}\label{e6}
\frac{\int_{B(y,\rho)}|\nabla u|^2}{\int_\Omega|\nabla u|^2}\le C\left(\frac{\int_{B(x,\rho)}|\nabla u|^2}{\int_\Omega|\nabla u|^2}\right)^{s}
\end{equation}
for all $y\in\Omega_{2\rho}$, where $0<s<1$ and $C$ depends on $\lambda_0,M_0,n,r_0,K_0,\rho,\Omega_\pm$.

In view of \eqref{e6} and covering $\Omega_{3\rho}$ with balls of radius $\rho$, we have that
\begin{equation}\label{e66}
\frac{\int_{\Omega_{3\rho}}|\nabla u|^2}{\int_\Omega|\nabla u|^2}\le C\left(\frac{\int_{B(x,\rho)}|\nabla u|^2}{\int_\Omega|\nabla u|^2}\right)^{s}.
\end{equation}
Note that $u-\phi_0$ is the solution to \eqref{aeq} with Dirichlet boundary value $\phi-\phi_0$. By Corollary 1.3 in \cite{livo}, we have that $$\|\nabla u\|_{L^{\infty}(\Omega)}^{2}=\|\nabla(u-\phi_0)\|_{L^\infty(\Omega)}^2\leq C\|\phi-\phi_0\|_{C^{1,\alpha'}(\partial\Omega)}^{2}$$ with $0<\alpha'<\frac{\alpha}{(\alpha+1)n}$, which implies
\begin{equation}
\int_{\Omega\setminus\Omega_{3\rho}}|\nabla u|^{2}\le C|\Omega\setminus\Omega_{5\rho}|\|\phi-\phi_0\|_{C^{1,\alpha'}(\partial\Omega)}^{2}\leq C\rho\|\phi-\phi_0\|_{C^{1,\alpha'}(\partial\Omega)}^{2}.
\label{ee3}
\end{equation}
Here we have used $\left|\Omega\backslash\Omega_{5\rho}\right|\lesssim\rho$ proved in \cite{ar}. Using the Poincar\'e inequality, we have
\begin{equation*}
\|\phi-\phi_0\|_{H^{1/2}(\partial\Omega)}^{2}\le C\|u-\phi_0\|_{H^{1}(\Omega)}^{2}\leq C\|\nabla u\|_{L^{2}(\Omega)}^{2}.\label{eq:e4}
\end{equation*}
Combining this and \eqref{ee3}, we see that if $\rho$ is small enough
depending on $\Omega_\pm$, $d_0$, $\lambda_0$, $M_{0}$, $r_0$, $K_0$, $s_0$, $L_0$, $\alpha$, $\alpha'$, $\rho$, and $\|\phi-\phi_0\|_{C^{1,\alpha'}(\partial\Omega)}/\|\phi-\phi_0\|_{H^{1/2}(\partial\Omega)}$, then
\[
\frac{\|\nabla u\|_{L^{2}(\Omega_{3\rho})}^{2}}{\|\nabla u\|_{L^{2}(\Omega)}^{2}}\ge\frac{1}{2}.
\]
The proposition follows from this and \eqref{e66}. \eproof

We now have enough tools to derive bounds on $|D|$.
\begin{theorem}\label{sizeest1} 
Suppose that the assumptions of this section hold. 
\begin{description}
\item[(i)] If, moreover, there exists $h>0$ such that
\begin{equation}
|D_{h}|\ge\frac{1}{2}|D|\quad\text{(fatness condition)}.\label{fatness}
\end{equation}
Then there exist constants $K_{1},K_{2}>0$ depending only on $\Omega_\pm$, $d_0$, $h$, $\lambda_0$, $M_{0}$, $r_0$, $K_0$, $s_0$, $L_0$, $\alpha$, $\alpha'$, and $\|\phi-\phi_0\|_{C^{1,\alpha'}(\partial\Omega)}/\|\phi-\phi_0\|_{H^{1/2}(\partial\Omega)}$, such that
\begin{equation*}
K_{1}\left|\frac{W_{0}-W}{W_{0}}\right|\le|D|\le K_{2}\left|\frac{W_{0}-W}{W_{0}}\right|.\label{size1001}
\end{equation*}

\item[(ii)] For a general inclusion $D$ contained strictly in $\Omega_+$, we assume that there exists $d_1>0$ such that
\[
\mbox{\rm dist}(D,\partial\Omega_+)\ge d_1.
\]
Then there exist constants $K_{1},K'_{2}$, $p>1$, depending only on $\Omega_\pm$, $d_0$, $d_1$, $h$, $\lambda_0$, $M_{0}$, $r_0$, $K_0$, $s_0$, $L_0$, $\alpha$, $\alpha'$, and $\|\phi-\phi_0\|_{C^{1,\alpha'}(\partial\Omega)}/\|\phi-\phi_0\|_{H^{1/2}(\partial\Omega)}$, such that
\begin{equation}
K_{1}\left|\frac{W_{0}-W}{W_{0}}\right|\le|D|\le K'_{2}\left|\frac{W_{0}-W}{W_{0}}\right|^{\frac 1p}.\label{size101}
\end{equation}
\end{description}
\end{theorem}
\pf. The proof follows closely the arguments in \cite{ars} and \cite{nw13}. The lower bound can be obtained by basic estimates. Let $c=\frac{1}{\left|\Omega_{d/4}\right|}\int_{\Omega_{d/4}}u$.
By the gradient estimate of \cite[Theorem 1.1]{livo}, the interior
estimate of \cite[Theorem 8.17]{gt} and the Poincar\'e inequality for
the domain $\Omega_{d/4}$, we have
\[
\|\nabla u\|_{L^{\infty}(\Omega_{d/2})}\leq C\|u-c\|_{L^{\infty}(\Omega_{d/3})}\leq C\|u-c\|_{L^{2}(\Omega_{d/4})}\leq C\|\nabla u\|_{L^{2}(\Omega)}.
\]
From this, the trivial estimate $\|\nabla u\|_{L^{2}(D)}^{2}\le C|D|\|\nabla u\|_{L^{\infty}(\Omega_{d/2})}^{2}$
and the second inequality of \eqref{ineq1}, the lower bound follows.

Next, we prove the upper bounds.

(i) Let $\rho=\frac{h}{4}$ and cover $D_{h}$
with internally nonoverlapping closed squares $\{Q_{k}\}_{k=1}^{J}$
of side length $2\rho$. It is clear that $Q_{k}\subset D$, hence
\begin{align*}
\int_{D}|\nabla u|^{2}dx & \ge\int_{\cup_{k=1}^{J}Q_{k}}|\nabla u|^{2}dx\ge\frac{|D_{h}|}{\rho^{2}}\min_{k}\int_{Q_{k}}|\nabla u|^{2}dx.\\
 & \geq\frac{C|D|}{\rho^{2}}\int_{\Omega}|\nabla u|^{2}dx.
\end{align*}
Here we have used Proposition \ref{lippro} and the fatness condition at the last inequality. The upper bound of $\left|D\right|$ follows
from this and the first inequality of \eqref{ineq1}.

(ii) To prove the upper bound without the fatness condition, we need the fact that $|\nabla u|^2$ is an $A_p$ weight which an easy consequence of the doubling condition for $\nabla u$.  It turns out when $D$ is strictly contained in $\Omega_+$ where the coefficient $A_+$ is Lipschitz. The well-known theorem guarantees that $|\nabla u|^2$ is an $A_p$ weight in $\Omega_+$ (see \cite{gl1} or \cite{ars}), i.e., for any $\bar r>0$, there exists $B>0$ and $p>1$ such that
\[
\left(\frac{1}{|B(a,r)|}\int_{B(a,r)}|\nabla u|^2\right)\left(\frac{1}{|B(a,r)|}\int_{B(a,r)}|\nabla u|^{-\frac{2}{p-1}}\right)^{p-1}\le B
\]
for any ball $B(a,r)\subset\Omega_{+,\bar r}$, where $B$ and $p$ depends on various constants listed in Proposition~\ref{lippro}. To derive the upper bound of \eqref{size101}, we choose $\bar r=d_1/2$ and follow exactly the same lines as in the proof of Theorem~2.2 \cite{ars}.

\eproof
\begin{remark}
We point out that part (i) of Theorem \ref{sizeest1} still holds if the assumption \eqref{interior} is replaced by
\[\mbox{dist}(D,\partial\Omega)\ge d_2>0.
\]
\end{remark}

\section*{Acknowledgements}
EF and SV were partially supported by GNAMPA - INdAM. EF was partially supported by the Research Project FIR 2013 Geometrical and
qualitative aspects of PDE's. JW was supported in part by MOST102-2115-M-002-009-MY3.

\end{document}